\pdfoutput=1

\documentclass[10pt]{article}

\makeatletter
\renewcommand{\maketitle}{\bgroup\setlength{\parindent}{0pt}
\begin{flushleft}
  \textbf{\@title}

  \@author
\end{flushleft}\egroup
} %%%%%
\def\footnoterule{\kern-3\p@
  \hrule \@width 2in \kern 2.6\p@} % the \hrule is .4pt high
\makeatother

\usepackage{preamble}
\usepackage{courier}
\usepackage[normalem]{ulem}
\usepackage[square,numbers]{natbib}

\newtheorem*{theorem1}{Theorem 1}
\newtheorem*{theorem2}{Theorem 2}
\newtheorem*{theorem3}{Theorem 3}

\title{{\Large \textbf{On non-abelian dp-minimal groups I: the torsion-free and distal cases}}}%\\{\large the torsion-free and distal cases}}}
\author{Atticus Stonestrom\ \ \ \ \ \ \ \ \ \ \ \ \ \ \ \ \ \  {\footnotesize (Email: \texttt{atticusstonestrom@yahoo.com})}}
\date{}

\begin{document}
\maketitle

\noindent\textbf{Abstract:} We give some results on dp-minimal groups. First we show that any torsion-free dp-minimal group is abelian; along the way show that any dp-minimal group admitting a principal f-generic type whose realizations are non-torsion is nilpotent-by-finite of class at most $2$. We then investigate the question of whether \textit{every} dp-minimal group $G$ is nilpotent-by-finite. There are naturally two cases: either (1) $G$ admits a distal f-generic type or (2) $G$ admits a generically stable f-generic type. In this paper we resolve case (1). This follows from a more general structural analysis, in which, assuming that $G$ admits a distal f-generic type, we show that the quotient of $G$ by its FC-center can be naturally equipped with the structure of a valued group; we then use this valuation structure to show that indeed $G$ is nilpotent-by-finite. Case (2) of the question will be resolved in an upcoming joint paper with Eran Alouf and Frank Wagner. \newline 

\noindent\textbf{Acknowledgements:} I would like to thank my PhD advisor Anand Pillay for all of his encouragement and feedback, and Caroline Terry for an inspirational discussion of the Gowers inverse question. \newline

\noindent\textbf{AI Declaration:} AI was not used for any part of this paper, neither in the mathematical content nor in the writing. \newline

\noindent\textbf{Notation:} Throughout $T$ will be a complete theory with a monster model $\mathfrak{C}$. For a parameter set $C\subseteq\mathfrak{C}$ and a type-definable set $X$, we use $S_X(C)$ to denote the space of all complete types over $C$ that are concentrated on $C$. For tuples $a,b$, and a small parameter set $C\subset\mathfrak{C}$, we write $a\equiv_C b$ to mean $\tp(a/C)=\tp(b/C)$.

Suppose $H$ is a group. We use the `left-normed' convention for commutators, so that $[a,b]=a^{-1}b^{-1}ab$ and $[a_1,\dots,a_{n+1}]=[[a_1,\dots,a_n],a_{n+1}]$. Conjugation is for us a right action, with $a^b=b^{-1}ab$ for all $a,b\in H$, and as usual we write $a^{-b}$ for $(a^{-1})^b$ and $a^{ib}=(a^i)^b$ for $i\in\mathbb{Z}$ and $a,b\in H$. If $H$ is a group we use $\sim_H$ to denote the relation of being conjugate in $H$. We use \uline{one piece of non-standard notation}: given sets $A_1,\dots,A_n$, we write $[A_1,\dots,A_n]$ to denote the set of commutators $\{[a_1,\dots,a_n]:a_i\in A_i\}$.

%We let $\mathrm{FC}(G)$ denote the FC-center of $G$, ie the subgroup of elements of $G$ whose $G$-conjugacy class is finite.

\section{Introduction}
`Dp-minimality', first introduced in \cite{shelah_strong_1} and \cite{shelah_strong_2} and then later isolated as a notion by Onshuus and Usvyatsov in \cite{onshuus_usvyatsov}, is a kind of abstract `one-dimensionality' condition which has proved a very fruitful strengthening of the NIP, both from the perspective of `pure' model theory and from the perspective of its algebraic implications; particularly notable examples of the latter are the classification of dp-minimal fields \cite{johnson} and the classification of dp-minimal integral domains \cite{delbee_halevi_johnson} (see also \cite{guingona},\cite{johnson_topology},\cite{jahnke_simon_walsberg},\cite{delbee_halevi}). Here we are concerned with dp-minimal groups. The typical examples, including for instance superstable groups of $U$-rank $1$, $(\mathbb{Z}/4\mathbb{Z})^\omega$, and the multiplicative and additive groups of $\mathbb{R},\mathbb{Q}_p,$ and algebraically closed valued fields, are all abelian-by-finite, but there exist dp-minimal groups which are not abelian-by-finite \cite{simonetta_example}, and it is an open question, asked in \cite{chernikov_pillay_simon} and \cite{kaplan_levi_simon}, whether every dp-minimal group is nilpotent-by-finite. This is known to be the case under several stronger hypotheses. For instance, C-minimal valued groups are nilpotent-by-finite \cite{simonetta_nilpotent}, $\omega$-categorical dp-minimal groups are abelian-by-finite \cite{dobrowolski_wagner}, profinite dp-minimal groups are abelian-by-finite \cite{profinite}, and left-ordered dp-minimal groups are abelian \cite{dobrowolski_goodrick}. In this paper we will give some more results. Let $G$ be a dp-minimal group. The first two results deal with consequences of torsion-free behavior:
\begin{theorem1} If $G$ admits a global bi-f-generic type concentrated on $G^{00}$ and whose realizations are torsion-free, then $G$ is nilpotent-by-finite of class at most $2$.
\end{theorem1}
\begin{theorem2}If $G$ is torsion-free then $G$ is abelian.
\end{theorem2}
These are Theorem \ref{gen_tf_implies_nil_2} and Theorem \ref{tf_implies_abelian} respectively. We will then turn to the more general question of whether dp-minimal groups are nilpotent-by-finite. Let $G$ be a dp-minimal group. In \cite{stonestrom} I proved that $G$ is definably amenable. By a result of Simon on dp-minimal invariant types, the question thus naturally splits into two cases: either (1) $G$ admits a distal f-generic type, or (2) $G$ admits a generically stable f-generic type. (See Section \ref{invariant types section} for definitions.) The third result of this paper is to resolve case (1):
\begin{theorem3}If $G$ admits a distal f-generic type, then $G$ is nilpotent-by-finite.
\end{theorem3}We will deduce Theorem 3 from a more general analysis, in which, assuming $G$ admits a distal f-generic type, we construct a valued group structure on the quotient of $G$ by its FC-center in such a way that the structure inherited by $G^{00}/Z(G^{00})$ is `plain' (see Section \ref{valued groups section}); this is the content of Section \ref{distal section}.

In fact we will show Theorems 1, 2, and 3 for the more general case when $G$ is a type-definable dp-minimal subgroup of an NIP group. This is for use in an upcoming joint paper with Eran Alouf and Frank Wagner which will resolve case (2), namely which will show that a dp-minimal group with a generically stable f-generic type is nilpotent-by-finite. Together with Theorem 3 it will follow that indeed every dp-minimal group is nilpotent-by-finite.

\section{Preliminaries}
\subsection{Groups}
Let $(H,\cdot)$ be an abstract group.

\subsubsection{Commutators}
Throughout this paper we will use the \uline{non-standard notation} that, for subsets $A_1,\dots,A_n\subseteq H$, $[A_1,\dots,A_n]$ denotes the set of commutators $\{[a_1,\dots,a_n]:a_i\in A_i\}$; when $A_i=\{a\}$ is a singleton we will just write $[A_1,\dots,A_{i-1},a,A_{i+1},\dots,A_n]$. The following is standard:
 \begin{lemma}\label{commutator_of_normal_subsets}
 Suppose $A_1,\dots,A_n\subseteq H$ are each fixed setwise under conjugation by elements of $H$. Then, for any $k\in\{0\}\cup[n]$, the set of $g\in H$ such that $[A_1,\dots,A_{k},g,A_{k+1},\dots,A_n]=\{1\}$ is a normal subgroup of $H$.
 \end{lemma}
 \begin{proof}
 Noting that $[A_1,\dots,A_{k},g]=[g^{-1},[A_1,\dots,A_{k}]^g]$, and that $[A_1,\dots,A_k]$ is fixed setwise under conjugation by elements of $H$, it suffices to show the claim when $k=0$, ie to show that the set of $g\in H$ with $[g,A_1,\dots,A_n]=\{1\}$ is a normal subgroup of $H$. This follows by induction on $n$; in the case $n=1$ the group is just $C_H(A_1)$.
 
For the inductive step, let $U$ denote the set of $u\in H$ with $[u,A_2,\dots,A_n]=\{1\}$ and let $V$ denote the set of $v\in H$ with $[v,A_1,\dots,A_n]=\{1\}$. Since $[v^h,A_1,\dots,A_n]=[v,A_1,\dots,A_n]^h$, $V$ is fixed setwise under conjugation by elements of $H$. Now, by the inductive hypothesis, $U$ is a normal subgroup of $H$.  If $g\in V$, then $[g,A_1]\subseteq U$, whence $[g^{-1},A_1]=[g,A_1]^{-g^{-1}}\subseteq U$, so that $g^{-1}\in V$. Finally, if $g,h\in V$, then $[g,A_1],[h,A_1]\subseteq U$, whence $[g,A_1]^h[h,A_1]\subseteq U$; since $[gh,A_1]\subseteq [g,A_1]^h[h,A_1]$ we thus have $gh\in V$, as desired.
 \end{proof}

 We also recall the Hall-Witt identity:

 \begin{fact} For all $a,b,c\in H$, $[c,b,a^c][a,c,b^a][b,a,c^b]=1$.
 \end{fact}
 
\subsubsection{Torsion-free groups}
We will need some facts about torsion-free groups; see for example 5.2.19 and 14.5.9 in \cite{robinson}.

\begin{fact}\label{tf_nilpotent_mod_center}
    Suppose $H$ is nilpotent and torsion-free. Then $H/Z(H)$ is torsion-free.
\end{fact}

\begin{fact}\label{tf_central_by_finite}
    If $H$ is torsion-free and $C_H(g)$ has finite index in $H$ for all $g\in H$, then $H$ is abelian.
\end{fact}

The following is an easy computation which we put here for convenience.
\begin{lemma}\label{conjugate_with_power_in_ab_by_torsion}
Suppose $a\in H$ is a non-torsion element such that, for some $i,j\in\mathbb{Z}$ and $g\in H$, we have $a^{ig}=a^j$. Suppose $H$ is abelian-by-torsion, or more generally that there is $n>0$ such that $[a^n,g^n]=1$. Then $|i|=|j|$.
\end{lemma}
\begin{proof}
    For notational convenience, let $f:H\to H$ denote the conjugation-by-$g$ map, so that $f(a^i)=a^j$ and $f^n(a^n)=a^n$. Let $b=a^n$; then $b$ is still non-torsion, and we have $f(b)^i=f(a^i)^n=a^{jn}=b^j$. By induction on $k$ it follows that $f^k(b)^{i^k}=b^{j^k}$ for all $k\geqslant 1$, since $$f^{k+1}(b)^{i^{k+1}}=f(f^k(b)^{i^{k}})^i=f(b^{j^k})^i=(f(b)^i)^{j^k}=b^{j^{k+1}}.$$ On the other hand, we have $b=f^n(b)$, and hence $b^{i^n}=f^n(b)^{i^n}=b^{j^n}$. Since $b$ is non-torsion, thus $i^n=j^n$, so that indeed $|i|=|j|$.
\end{proof}

\subsubsection{Valued groups}\label{valued groups section}
Let $(\mathcal{P},\leqslant)$ be a totally ordered set with a greatest element. A $\mathcal{P}$-valued group structure on $H$ is a a surjective map $v:H\to\mathcal{P}$ satisfying the following:
\begin{enumerate}
\item $v(a)=\max\mathcal{P}$ if and only if $a=1$.
\item $v(a^g)<v(b^g)$ for all $g\in H$ and $a,b\in H$ with $v(a)<v(b)$.
%\item $v(ab^{-1})\geqslant\min\{v(a),v(b)\}$.
\item $v(a)=v(a^{-1})$ for all $a\in H$
\item $v(ab)\geqslant\min\{v(a),v(b)\}$ for all $a,b\in H$.
\end{enumerate}
Following \cite{simonetta_nilpotent}, we will say that $(H,v)$ is `plain' if we additionally have $v([a,b])>\max\{v(a),v(b)\}$ for all $a,b\in H\setminus\{1\}$. The following is proved in \cite{simonetta_nilpotent}:
\begin{fact}\label{plain_implies_loc_nil}
Suppose $(H,v)$ is a plain valued group of finite exponent. Then $H$ is locally nilpotent.
\end{fact}

\subsection{Model theory}
Throughout this section, \uline{assume that $T$ is NIP and that $(\tilde{G},\cdot)$ is a definable group of $T$}. If $G$ is a type-definable subgroup of $\tilde{G}$, by a `relatively definable' subset of $G$ we mean a type-definable set of form $X\cap G$ for some definable set $X\subseteq\tilde{G}$. As usual, if $G$ is a type-definable subgroup of $\tilde{G}$, then whenever it does not cause confusion we will identify $G$ with the abstract group of its $\mathfrak{C}$-points $G(\mathfrak{C})$. The following is an easy consequence of compactness:

\begin{fact}\label{bounded index rel definable implies finite index}
    Let $G$ be a type-definable subgroup of $\tilde{G}$ and let $H\leqslant G$ be a relatively definable subgroup of $G$. If $H$ has bounded index in $G$, then $H$ has finite index in $G$.
\end{fact}

\subsubsection{NIP groups}
We will use several general facts about NIP groups. First recall the Baldwin-Saxl theorem, proved in \cite{baldwin_saxl}: for any $L$-formula $\phi(x,y)$, there is a constant $k$ such that, for all $c_1,\dots,c_n\in\mathfrak{C}^y$ such that each $\phi(x,c_i)$ defines a subgroup of $\tilde{G}$, there is a subset $s\subseteq[n]$ of size at most $k$ with $\bigcap_{i\in s}\phi(x,c_i)=\bigcap_{i\in [n]}\phi(x,c_i)$. We use the following consequence of this:
\begin{fact}\label{baldwin_saxl_centralizers}
    There is a constant $n$ such that, for any finite subset $B\subset \tilde{G}(\mathfrak{C})$, there is $S\subseteq B$ of size at most $n$ with $C_{\tilde{G}}(S)\subseteq C_{\tilde{G}}(b)$ for all $b\in B$. In particular, if $(u_i)_{i\in[n]}+(a_i)_{i\in\omega}+(v_i)_{i\in[n]}$ is indiscernible and $W=\{u_i,v_i:i\in[n]\}$, then $C_{\tilde{G}}(W)\subseteq C_{\tilde{G}}(a_i)$ for all $i\in \omega$.
\end{fact}

We will also make extensive use of the `existence' of $G^{00}$ \cite{shelah_g00} and $G^{\infty}$ \cite{gismatullin}:

\begin{fact}\label{g00_exists}
    Let $G\leqslant\tilde{G}$ be a subgroup type-definable over a small set $A$. Then there is a bounded-index subgroup $G^{00}\leqslant G$ (resp. $G^\infty\leqslant G$) type-definable (resp. invariant) over $A$ and contained in any bounded-index subgroup of $G$ type-definable (resp. invariant) over a small set of parameters.
\end{fact}

We will also need some of the theory of definably amenable NIP groups, as developed in \cite{hrushovski_peterzil_pillay}, \cite{hrushovski_pillay}, \cite{chernikov_simon}, and \cite{stonestrom}. \uline{Let $G(x)$ be a type-definable subgroup of $\tilde{G}$.} We say $G$ is `definably amenable' if there is a $G$-invariant Keisler measure concentrated on $G$, ie a Keisler measure $\mu$ with $\mu(X)=1$ for every definable supserset $X\supseteq G$ and such that $\mu(gX)=\mu(X)$ for every definable set $X\subseteq\tilde{G}$ and every $g\in G$. Given a partial type $\pi(x)$ concentrated on $\tilde{G}$, we say (1) $\pi(x)$ is (left) `f-generic for $G$' if, for any small model $M$ over which $\pi$ and $G$ are type-definable, no left translate of the partial type $\pi(x)\wedge G(x)$ by an element of $G$ forks over $M$, and (2) $\pi(x)$ (left) `$G$-divides' if there are definable sets $X,Y\subseteq\tilde{G}$ with $\pi(x)\vdash x\in X$ and $G(x)\vdash x\in Y$, and a sequence of elements $(g_i:i\in\omega)$ from $G$, such that, for some $k\in\omega$, the family of translates $(g_i(X\cap Y):i\in\omega)$ is $k$-inconsistent. Now:
\begin{fact}\label{f-generic facts}
    \begin{enumerate}
        \item A partial type $\pi(x)$ concentrated on $\tilde{G}$ is f-generic for $G$ if and only if it does not $G$-divide.
        \item $G$ is definably amenable if and only if there is a complete global type f-generic for $G$ if and only if the definable subsets of $\tilde{G}$ not f-generic for $G$ form an ideal.
        \item If $G$ is definably amenable then, for any small model $M\prec\mathfrak{C}$, there is a global $M$-invariant type f-generic for $G$.
        \item\label{f-gen_iff_g00_inv} A complete global type concentrated on $G$ is f-generic for $G$ if and only if it is $G^{00}$-invariant. In particular:
        \item\label{g00=g000}Suppose $p\in S_G(\mathfrak{C})$ is f-generic for $G$. Then, for any small $D\subset\mathfrak{C}$, every element of $G^{00}(\mathfrak{C})$ is a difference of two realizations of $p|_D$, ie of form $ab^{-1}$ for some $a,b\models p|_D$.
    \end{enumerate}
\end{fact} 

A consequence of Fact \ref{f-generic facts}(\ref{g00=g000}) is the following:

\begin{fact}\label{f-generic invariant subgroup contains G00}
    Suppose $p\in S_G(\mathfrak{C})$ is f-generic for $G$, and suppose that $H$ is a subgroup of $G$, invariant over some small parameter set $D$, such that some (equivalently any) realization of $p|_D$ lies in $H$. Then $G^{00}\leqslant H$.
\end{fact}

Note that, if $p(x)\in S_{G}(\mathfrak{C})$ is invariant over a small model and f-generic for $G$, then the pushforward of $p_x\otimes p_y$ under the map $(x,y)\mapsto xy^{-1}$ is (i) concentrated on $G^{00}$, and (ii) invariant under both left-translation and right-translation by all elements of $G^{00}(\mathfrak{C})$. Thus, defining right-f-genericity and bi-f-genericity in the natural way, one has:

\begin{fact}\label{bi-f-gens_exist}
If $G$ is definably amenable, then there is a complete global type invariant over a small model, concentrated on $G^{00}$, and bi-f-generic for $G$.
\end{fact}

\subsubsection{Dp-minimality}
There are numerous equivalent characterizations of dp-minimality; see for example \cite{dolich_goodrick_lippel} for an introduction. We use one coming from \cite{usvyatsov_gs}; recall that sequences $I_1,\dots,I_n$ are said to be `mutually indiscernible' over a parameter set $B$ if, for each $k\in[n]$, the sequence $I_k$ is indiscernible over $(B,I_{\neq k})$.

\begin{definition}
    Given a natural number $n$, a partial type $\pi(x)$ over a small parameter set is said to have dp-rank at most $n$ if the following holds for some (equivalently any) small parameter set $B$ over which $\pi$ is type-defined: for any infinite sequences $I_1,\dots,I_{n+1}$ mutually indiscernible over $B$, and any realization $a\models\pi(x)$, one of $I_1,\dots,I_{n+1}$ is indiscernible over $(a,B)$. A global type $p(x)$ has dp-rank at most $n$ if there is a small parameter set $B$ such that $p(x)|_B$ has dp-rank at most $n$. Finally, a partial type or global type is said to be dp-minimal if it has dp-rank at most $1$.
\end{definition}

Then a partial type $\pi(x)$ is not dp-minimal if and only if there are $L(\mathfrak{C})$-formulas $\phi(x,y)$ and $\psi(x,z)$ and tuples $b_i\in\mathfrak{C}^y$ and $c_j\in\mathfrak{C}^z$ such that, for all $i,j\in\omega$, the partial type $$\pi(x)\wedge\phi(x,b_i)\wedge\psi(x,c_j)\wedge\bigwedge_{k\neq i}\neg\phi(x,b_k)\wedge\bigwedge_{l\neq j}\neg\psi(x,c_l)$$ is consistent.

Now, \uline{let $G$ be a type-definable subgroup of $\tilde{G}$}. We say that the $G$ is dp-minimal if the partial type $G(x)$ is dp-minimal. The following was noted in \cite{simon_ordered} for definable subgroups of definable dp-minimal groups, and the analogue is true for invariant subgroups:

\begin{fact}\label{inv_subgps_are_commensurable}
Suppose $G$ is dp-minimal, and that $A,B\leqslant G$ are subgroups invariant over some small model. Then $A\cap B$ has bounded index in one of $A$ and $B$.
\end{fact}
\begin{proof}
Let $M$ be a small model over which $G$ is type-definable and $A,B$ are invariant. Suppose for contradiction that $A\cap B$ has unbounded index in both $A$ and $B$. Then there are sequences $I=(a_i)_{i\in\omega}$ in $A$ and $J=(b_j)_{j\in\omega}$ in $B$, each indiscernible over $M$, and such that $a_iB\cap a_jB=Ab_i\cap Ab_j=\varnothing$ for all $i\neq j$. By Ramsey and $M$-invariance of $A,B$ we may assume that $I$ and $J$ are mutually indiscernible over $M$, and then by dp-minimality one of $I$ and $J$ must be indiscernible over $(M,a_0b_0)$. In the former case, there is an automorphism fixing $(M,a_0b_0)$ and taking $a_0$ to $a_1$, and since $a_0^{-1}(a_0b_0)\in B$ it follows by invariance that $a_1^{-1}(a_0b_0)\in B$. But then $a_1^{-1}a_0\in B$, a contradiction. A symmetric argument, using that $(a_0b_0)b_0^{-1}\in A$, gives the result in the latter case.
\end{proof}

We will also use the following, from \cite{stonestrom} and \cite{simon_ordered} respectively. In both cases the results are written for definable dp-minimal groups, but the same arguments work for type-definable ones; the proof of Fact \ref{dp-min_implies_def_am} works as written in \cite{stonestrom} for the type-definable case, but the proof of Fact \ref{dp-min_implies_ab_by_tor} from \cite{simon_ordered} needs a small amount of rewriting to work for the type-definable case, which we give for completeness.

\begin{fact}\label{dp-min_implies_def_am}
Any dp-minimal group is definably amenable.
\end{fact}

\begin{fact}\label{dp-min_implies_ab_by_tor}
    Suppose $G$ is dp-minimal. Then there is a definable superset $X\supseteq G$ and some $n>0$ such that, for all $a,b\in X$, one of $[a^n,b]=1$ and $[a,b^n]=1$ holds. In particular, $G$ is abelian-by-torsion.
\end{fact}
\begin{proof}
    If $A,B$ are relatively definable subgroups of $G$, then by Fact \ref{inv_subgps_are_commensurable} $A\cap B$ has bounded index in one of $A$ and $B$, and then by Fact \ref{bounded index rel definable implies finite index} $A\cap B$ has finite index in one of $A$ and $B$. Now, for any $x\in G$, $C_{G}(x)=C_{\tilde{G}}(x)\cap G$ is a relatively definable subgroup of $G$. So, for any $x,y\in G$, $C_G(x)\cap C_G(y)$ has finite index in one of $C_G(x)$ and $C_G(y)$, so in particular there must be some $k$ such that $y^k\in C_G(x)$ or $x^k\in C_G(y)$. So the partial type $$G(x)\wedge G(y)\wedge\{[x^k,y]\neq 1,[x,y^k]\neq 1:k\in\omega\}$$ is inconsistent, and so by compactness there is a definable superset $X\supseteq G$ and some $k\in\omega$ such that, for all $a,b\in X$, there is $l\leqslant k$ such that one of $[a^{l},b]=1$ and $[a,b^{l}]=1$ holds. Then, letting $n=k!$, for all $a,b\in X$, one of $[a^{n},b]=1$ and $[a,b^{n}]=1$ holds, yielding the first part of the statement.

    For the second part of the statement, let $Y=\{x^n:x\in X\}$, a definable set. Then any two elements of $Y$ commute with each other, so $\tilde{A}:=C_{\tilde{G}}(C_{\tilde{G}}(Y))$ is a definable abelian subgroup of $\tilde{G}$ containing $Y$. So now $A:=\tilde{A}\cap G$ is a relatively definable abelian subgroup of $G$, and since $G\subseteq X$ we have $g^n\in Y\subseteq\tilde{A}$ for all $g\in G$, so $g^n\in A$ for all $g\in G$. 

    When we say that $G$ is abelian-by-torsion we are as usual identifying the type-definable group $G$ with the abstract group of its $\mathfrak{C}$-points $G(\mathfrak{C})$. Now $A(\mathfrak{C})$ is an abelian subgroup of $G(\mathfrak{C})$ containing $g^n$ for all $g\in G(\mathfrak{C})$. In particular, the `normal core' $\bigcap_{g\in G(\mathfrak{C})}A(\mathfrak{C})^g$ still contains $g^n$ for all $g\in G(\mathfrak{C})$, and it is a normal abelian subgroup of $G(\mathfrak{C})$, so $G(\mathfrak{C})$ is abelian-by-(exponent $\leqslant n$).
\end{proof}

A fact we will not need but which is occasionally useful for discussion is the following, from \cite{kaplan_levi_simon}; again the result there is written only for definable dp-minimal groups, but the argument works for type-definable ones as well. Since we do not rely on this fact for our results we do not include include details.

\begin{fact}\label{f-gen_equals_strong_f-gen}
If $G$ is dp-minimal and $q(x)\in S_G(\mathfrak{C})$ is f-generic for $G$, then $q$ is $M$-invariant for every small $M\prec\mathfrak{C}$.
\end{fact}

\subsubsection{Invariant types in NIP theories}\label{invariant types section}
In this section continue to assume that $T$ is NIP. We will need several facts about invariant types under this assumption. Recall from \cite{shelah_strong_1} that, if $M\prec\mathfrak{C}$ is a small model and $q(x)\in S(\mathfrak{C})$ is a global type, then $q(x)$ does not fork over $M$ if and only if it is $M$-invariant. Thus invariant types are plentiful in the NIP context.

Recall the notion of a `generically stable' type, introduced in \cite{shelah_gs}, then investigated in \cite{hrushovski_pillay} and in \cite{usvyatsov_gs} (mostly with an NIP hypothesis) and in \cite{pillay_tanovic} (without an NIP hypothesis). A global type $p(x)\in S(\mathfrak{C})$ invariant over some small set is said to be `generically stable' if for some (any) small $B\subset\mathfrak{C}$ over which it is invariant, we have: for any linear order $I$, any Morley sequence $(a_i)_{i\in I}\models p^{\otimes I}|_B$, and any $L(\mathfrak{C})$-formula $\phi(x,c)$, the set $\{i:a_i\models\phi(x,c)\}$ is either finite or co-finite in $I$. There are many equivalent properties under the NIP hypothesis; for example $p$ is generically stable if and only if it is definable and finitely satisfiable in some small model.

The other class of invariant type we will use is that of distal types, implicitly introduced in \cite{simon_distality},\footnote[1]{In an NIP theory, distal types are to distal indiscernible sequences as generically stable types are to totally indiscernible sequences.} but first explicitly introduced in \cite{simon_book}. For the definition, recall that, for a parameter set $B\subset\mathfrak{C}$ and disjoint tuples of variables $x,y$, two types $p(x)\in S_x(B)$ and $q(y)\in S_y(B)$ are said to be `weakly orthogonal', denoted $p(x)\perp^w q(y)$, if $p(x)\wedge q(y)$ determines a unique complete type in $S_{(x,y)}(B)$. Now a $B$-invariant global type $p(x)\in S(\mathfrak{C})$ is said to be `distal over $B$' if the following holds: for any tuple $c$, if $I\models p^{\otimes\omega}|_{(B,c)}$ then $\tp(c/B,I)\perp^w p|_{(B,I)}$. It is proved in \cite{simon_book} that, if $p(x)$ is a global type invariant over sets $B,C\subset\mathfrak{C}$, then $p(x)$ is distal over $B$ if and only if it is distal over $C$. So one may speak of a `distal type' without ambiguity.

Distal types and generically stable types lie at precisely opposite ends of the NIP spectrum, respectively capturing `purely order-like' and `purely stable-like' behavior in NIP theories. Simon showed that, in the dp-minimal case, every 1-type falls into one of those two classes; this is again morally proved in \cite{simon_distality}, but for the explicit formulation we use see Proposition 9.13 in \cite{simon_book}:
\begin{fact}\label{dp-min_implies_gs_or_distal}
Suppose $B\subset\mathfrak{C}$ is a small set, and that $p(x)\in S(\mathfrak{C})$ is a $B$-invariant global type such that $p(x)|_B$ is dp-minimal. Then $p(x)$ is either generically stable or distal.
\end{fact}

Now, still under the global NIP assumption, suppose that $G$ is a type-definable group, and suppose there is a global generically stable type $p\in S_G(\mathfrak{C})$ f-generic for $G$. Then every global f-generic type for $G$ is a translate of $p$. To see this, suppose $q\in S_G(\mathfrak{C})$ is an f-generic type for $G$. Now $p,q$ respectively concentrate on some cosets $aG^{00},bG^{00}$ of $G^{00}$, and then $a^{-1}p,b^{-1}q$ are f-generic types of $G$ concentrated on $G^{00}$. In particular by Fact \ref{f-generic facts}(\ref{f-gen_iff_g00_inv}) both $a^{-1}p$ and $b^{-1}q$ are $G^{00}$-invariant. Since $p$ is generically stable, also $a^{-1}p$ is generically stable, and by Lemma 1 in \cite{pillay_tanovic} if follows that $a^{-1}p=b^{-1}q$, so indeed $q=ba^{-1}p$. So combining this observation with Fact \ref{dp-min_implies_gs_or_distal} and Fact \ref{f-gen_equals_strong_f-gen} now yields the following:

\begin{corollary}\label{dp-min_implies_gs_or_distal_group_version}
Suppose $G$ is a type-definable dp-minimal group. Then either every global type f-generic for $G$ is generically stable, or every global type f-generic for $G$ is distal.
\end{corollary}

\section{Nilpotence in a Morley sequence}
In this section \uline{assume that $T$ is NIP, that $\tilde{G}$ is a definable group of $T$, and that $G$ is a type-definable subgroup of $\tilde{G}$}. We will look at consequences of certain commutator relations on a Morley sequence of a global invariant types. Despite not being proved explicitly anywhere in the literature the material in Section \ref{generic nilpotence section} is very standard. On the other hand I hope that the material of Section \ref{local nilpotence section} is more novel/interesting.

\subsection{Generic nilpotence}\label{generic nilpotence section}

It is well known that, in a stable group, to check virtual nilpotence of class at most $(c-1)$ it suffices to check that an independent $c$-tuple of principal generics satisfies $[x_1,\dots,x_c]=1$; see eg \cite{wagner}. We point out here an analogous fact in the NIP setting, which follows quickly from Fact \ref{g00_exists}. The proof is standard; compare also de Aldama's proof of the existence of nilpotent envelopes in NIP groups \cite{de_aldama}.

\begin{lemma}\label{generic_nil_is_nil}
    Let $p(x)\in S_{G^{00}}(\mathfrak{C})$ be f-generic for $G$ and invariant over some small model. Suppose $p^{\otimes c}(x_1,\dots,x_c)\vdash [x_1,\dots,x_c]=1$. Then there is a definable subgroup $O\leqslant \tilde{G}$ such that $O$ is nilpotent of class at most $(c-1)$ and such that $O\cap G$ has finite-index in $G$.
\end{lemma}
\begin{proof}
    By induction on $c$. In the case $c=1$, ie when $p(x)\vdash x=1$, then $G$ must be finite and the result follows taking $O=\{1\}$. Now suppose the result holds for $c$, and that $p(x)\in S_{G^{00}}(\mathfrak{C})$ is an invariant type f-generic for $G$ with $p^{\otimes {(c+1)}}(x_1,\dots,x_{c+1})\vdash[x_1,\dots,x_{c+1}]=1$. Let $q(y)$ be the pushforward of $p^{\otimes c}(x_1,\dots,x_c)$ under the map $(x_1,\dots,x_c)\mapsto [x_1,\dots,x_c]$; then $p(x)\otimes q(y)\vdash [x,y]=1$. 

    Let $M$ be a small model over which $\tilde{G}$ is definable, $G$ and hence $G^{00}$ are type-definable, and $p$ and hence $q$ are invariant. Let $b\models q|_M$. Now $p(x)\vdash [x,b]=1$. Hence $C_G(b)=C_{\tilde{G}}(b)\cap G$ is a relatively-definable subgroup of $G$ concentrated on by $p$; since $p$ is f-generic for $G$ it follows from Fact \ref{f-generic invariant subgroup contains G00} that $G^{00}\leqslant C_G(b)$. In particular $G^{00}(x)\vdash x\in C_{\tilde{G}}(b)$, so by compactness there is an $M$-definable superset $D\supseteq G^{00}$ such that $D\subseteq C_{\tilde{G}}(b)$. Then $b\in C_{\tilde{G}}(D)$, and since $D$ is $M$-definable and $b\models q|_M$ thus $q(y)\vdash y\in C_{\tilde{G}}(D)$.

    Now, let $H=C_{\tilde{G}}(C_{\tilde{G}}(D))$. Then $H\supseteq D\supseteq G^{00}$. So $H\cap G$ is a relatively-definable subgroup of $G$ containing $G^{00}$, hence by Fact \ref{bounded index rel definable implies finite index} $H\cap G$ has finite-index in $G$; call this ($\ast$). Since $G^{00}\leqslant H$, both $p$ and $q$ are concentrated on $H$. Also, since $q$ is concentrated on $C_{\tilde{G}}(D)$, in fact $q$ is concentrated on $Z(H)$. Call this ($\ast\ast$).

    Working in $T^{\mathrm{eq}}$ we may identify $H/Z(H)$ with a definable group. Let $\bar{p}$ be the pushforward of $p$ under the quotient map $H\to H/Z(H)$; this is well-defined since $p(x)\vdash x\in H$. Let $\overline{G\cap H}$ be the pushforward of $G\cap H$ under the quotient map; then it is easy to check that $\overline{G\cap H}$ is a type-definable subgroup of $H/Z(H)$ and that $\bar{p}$ is an f-generic type for $\overline{G\cap H}$ concentrated on $\overline{G\cap H}{}^{00}$. By ($\ast\ast$) we have $\bar{p}{}^{\otimes c}(x_1,\dots,x_c)\vdash [x_1,\dots,x_c]=1$, so by the inductive hypothesis there is a definable subgroup $P\leqslant H/Z(H)$ such that $P$ is nilpotent of class at most $(c-1)$ and such that $P\cap\overline{G\cap H}$ has finite-index in $\overline{G\cap H}$. Let $O\leqslant H$ be the preimage of $P$ under the quotient map $H\to H/Z(H)$. Now $Z(H)\leqslant Z(O)$, so $O/Z(O)$ is a quotient of $O/Z(H)=P$ and hence $O$ is nilpotent of class at most $c$. Since $P\cap\overline{G\cap H}$ has finite-index in $\overline{G\cap H}$, $O\cap ((G\cap H)Z(H))$ has finite-index in $(G\cap H)Z(H)$. Taking intersections with $G\cap H$ gives that $O\cap G\cap H=O\cap G$ has finite-index in $G\cap H$. By ($\ast$), thus $O\cap G$ has finite-index in $G$.
\end{proof}

\subsection{Local nilpotence in dp-finite groups}\label{local nilpotence section}
We will only need the following result on dp-finite groups in the dp-minimal case, but I include the more general result in the hope that it may be of broader interest or use. Compare also Dobrowolski and Goodrick's chain condition for finite burden groups from \cite{dobrowolski_goodrick}. 
\begin{lemma}\label{gen_n-nil_implies_gen_N_nil}
    Suppose that $G$ has dp-rank at most $n-1$. Let $M\prec\mathfrak{C}$ be any model, and suppose $p(x)\in S_G(\mathfrak{C})$ is an $M$-invariant type with the property that, for some $N$, the group $\langle a_1,\dots,a_n\rangle$ is nilpotent of class at most $(N-1)$ for some (any) $(a_1,\dots,a_n)\models p^{\otimes n}|_M$. Then there is some $\sigma\in S_N$ such that $p^{\otimes N}(x_1,\dots,x_N)\vdash [x_{\sigma(1)},\dots,x_{\sigma(N)}]=1$.
\end{lemma}
\begin{proof}
    Let $I_1+J_1+\dots+I_n+J_n$ be a Morley sequence of $p$ over $M$, with $I_k=(a_{(k,i)})_{i\in\omega}$ and $J_k=(b_{(k,j)})_{j\in\mathbb{Z}}$. Write $J=J_1+\dots+J_n$, ie $J=(b_{\alpha})_{\alpha\in [n]\times\mathbb{Z}}$, so that $J$ is a Morley sequence of $p$ of order type $[n]\times\mathbb{Z}$, ordered lexicographically. Now by the hypothesis of the lemma $\langle a_{(1,0)},a_{(2,0)},\dots,a_{(n,0)}\rangle$ is nilpotent of class at most $(N-1)$.

    We will prove by backwards induction on $k$ that, for every $k\in[N]\cup\{0\}$, and for all $c_1,\dots,c_k\in\{a_{(1,0)},a_{(2,0)},\dots,a_{(n,0)}\}$, there are pairwise distinct $\alpha_1,\dots,\alpha_{N-k}\in [n]\times\mathbb{Z}$ such that $[c_1,\dots,c_k,b_{\alpha_{N-k}},\dots,b_{\alpha_{1}}]=1$. At the end of the induction, when $k=0$, we will have the desired result. The base case $k=N$ follows from the assumption that $\langle a_{(1,0)},a_{(2,0)},\dots,a_{(n,0)}\rangle$ is nilpotent of class at $(N-1)$.

    Now suppose the result holds for $k+1$, and fix $c_1,\dots,c_k\in\{a_{(1,0)},a_{(2,0)},\dots,a_{(n,0)}\}$. $(I_1,\dots,I_n)$ are mutually indiscernible over $J$, so by the dp-rank assumption there is $\ell\in[n]$ such that $I_\ell$ is indiscernible over $(J,[c_1,\dots,c_k])$. By there inductive hypothesis applied with $c_{k+1}=a_{(\ell,0)}$, there are pairwise distinct $\alpha_1,\dots,\alpha_{N-k-1}\in [n]\times\mathbb{Z}$ such that $$[c_1,\dots,c_k,a_{(\ell,0)},b_{\alpha_{N-k-1}},\dots,b_{\alpha_1}]=1.$$ Since $I_\ell$ is indiscernible over $([c_1,\dots,c_k],b_{\alpha_1},\dots,a_{\alpha_{N-k-1}})$, in particular also $$[c_1,\dots,c_k,a_{(\ell,1)},b_{\alpha_{N-k-1}},\dots,b_{\alpha_1}]=1.$$ Now pick some $j\in\mathbb{Z}$ such that $j<j'$ for every $j'\in\mathbb{Z}$ such that $(\ell,j')$ appears among the $\alpha_1,\dots,\alpha_{N-k-1}$. By definition of the $c_i$ and indiscernibility of $I_1+J_1+\dots+I_n+J_n$, then $a_{(\ell,1)}$ and $b_{(\ell,j)}$ have the same type over $([c_1,\dots,c_k],b_{\alpha_1},\dots,a_{\alpha_{N-k-1}})$, so $$[c_1,\dots,c_k,b_{(\ell,j)},b_{\alpha_{N-k-1}},\dots,b_{\alpha_1}]=1,$$ and thus we may take $\alpha_{N-k}=(\ell,j)$.
\end{proof}

% \subsection{Generically 2-Engel}
% Suppose now that $G$ is dp-minimal.

\section{Torsion-free case}\label{torsion-free section}
In this section \uline{assume that $T$ is NIP, that $\tilde{G}$ is a definable group of $T$, and that $G$ is a type-definable dp-minimal subgroup of $\tilde{G}$}. Additionally, by adding constant symbols, for notational convenience \uline{assume that $G$ is type-definable without parameters}. We will show that, if $G$ has a global bi-f-generic type concentrated on $G^{00}$ whose realizations are non-torsion, then $G$ is nilpotent-by-finite of class at most $2$. Moreover, if $G$ is torsion-free, then $G$ is abelian.
 
\begin{lemma}\label{seqs_of_non-conjugates_commute}
Suppose that $I=(a_i)_{i\in\omega}$ is a sequence in $G$ such that (i) for some definable superset $X\supseteq G$, we have $a_i^g\neq a_j$ for all $g\in X$ and $i\neq j$, and (ii) $[a_i,a_j]=1$ for all $i,j$. Then, for every sequence $J=(b_j)_{j\in\omega}$ in $G$, there is some $i\neq j$ and $k\in\omega$ with $[a_k^{b_i},a_k^{b_j}]=1$.
\end{lemma}
\begin{proof}
    Suppose that $J$ is otherwise, so that $[a_k^{b_i},a_k^{b_j}]\neq 1$ for all $i\neq j$ and $k\in\omega$. By Ramsey and compactness, we may then assume that $I$ and $J$ are mutually indiscernible sequences. Now, let $n$ be given Fact \ref{baldwin_saxl_centralizers}. Pick $u_i,v_i$ so that $(u_i)_{i\in[n]}+I+(v_i)_{i\in[n]}$ is still $J$-indiscernible and let $W=\{u_i,v_i:i\in[n]\}$; then $I$ and $J$ are mutually indiscernible over $W$. Moreover, $D:=C_{\tilde{G}}(W)$ is $W$-definable, and by condition (ii) and Fact \ref{baldwin_saxl_centralizers} we have $a_k\in D\subseteq C_{\tilde{G}}(a_k)$ for all $k\in\omega$.
    
    Now, by condition (i), we have $a_i\sim a_{0}^{b_0}$ iff $i=0$, so $I$ is not indiscernible over $a_0^{b_0}$ and hence not over $(W,a_0^{b_0})$. Thus, by dp-minimality, $J$ must be indiscernible over $(W,a_0^{b_0})$. Since $a_0\in D$, we have $b_0\models a_0^{b_0}\in D^y$, and thus also $a_0^{b_0}\in D^{b_1}$. But $D\subseteq C_{\tilde{G}}(a_0)$, so that $a_0^{b_0}\in C_{\tilde{G}}(a_0)^{b_1}=C_{\tilde{G}}(a_0^{b_1})$, ie $[a_0^{b_0},a_0^{b_1}]=1$, a contradiction.
\end{proof}

\begin{lemma}\label{tf_commutes_with_ind_conjugates}
    Suppose that $a\in G(\mathfrak{C})$ is non-torsion, and that $J=(g_i)_{i\in\omega}$ is an $a$-indiscernible sequence in $G$. Then $[a^{g_0},a^{g_1}]=1$.
\end{lemma}
\begin{proof}
    By Fact \ref{dp-min_implies_ab_by_tor}, there is a definable superset $X\supseteq G$ and some $n>0$ such that, for all $a,b\in X$, we have $[a^n,b^n]=1$. In particular, since $a$ is non-torsion, by Lemma \ref{conjugate_with_power_in_ab_by_torsion} we have $a^{ig}\neq a^j$ for all $g\in X$ and $i\neq j\in\omega$. So any infinite sequence of positive powers of $a$ satisfies the hypotheses of Lemma \ref{seqs_of_non-conjugates_commute}.

    Let $s\subseteq\omega$ be the set of all $k\in\omega$ such that $[a^{kg_0},a^{kg_1}]\neq 1$; since $J$ is $a$-indiscernible, note that, if $k\in s$, then $[a^{kg_i},a^{kg_j}]\neq 1$ for all $i\neq j$. So applying Lemma \ref{seqs_of_non-conjugates_commute} to the sequences $I=(a^i)_{i\in s}$ and $J=(g_i)_{i\in\omega}$ it follows that $s$ must be finite. On the other hand, letting $b=a^{g_0}$ and $c=a^{g_1}$, we have $[b^k,c^k]=1$ for all $k\notin s$, and it follows that $[b,c]=1$. (For example, since $s$ is finite, we may find four pairwise coprime numbers $p,q,r,s\in\omega\setminus s$. Then $b^p$ and $b^q$ each lie in $C_{\tilde{G}}(c^{pq})$, so $b\in C_{\tilde{G}}(c^{pq})$, and likewise $b\in C_{\tilde{G}}(c^{rs})$. Now $c^{pq}$ and $c^{rs}$ each lie in $C_{\tilde{G}}(b)$, so $c\in C_{\tilde{G}}(b)$, as needed.)
\end{proof}

\begin{lemma}\label{tf_f-gen_is_2-engel}
    Let $M$ be a small model and suppose that $q(y)\in S_{G^{00}}(\mathfrak{C})$ is an $M$-invariant global type right f-generic for $G$. Then $q(y)\vdash [a,a^y]=1$ for any non-torsion element $a\in G(\mathfrak{C})$. Moreover, if $p(x)\in S_G(\mathfrak{C})$ is an $M$-invariant type whose realizations are non-torsion, then $p(x)\otimes q(y)\vdash [x,x^y]=1$. (And by the first part also $q(y)\otimes p(x)\vdash [x,x^y]=1$.)
\end{lemma}
\begin{proof}
    For the first part, let $a\in G(\mathfrak{C})$ be a non-torsion element and let $J=(g_i)_{i\in\omega}\models q^{\otimes\omega}|_{(M,a)}$. Then $J$ is $a$-indiscernible, so by Lemma \ref{tf_commutes_with_ind_conjugates} $[a^{g_0},a^{g_1}]=1$, hence $[a,a^{g_1g_0^{-1}}]=1$. But $q$ is right f-generic for $G$ and concentrated on $G^{00}$, so $g_1g_0^{-1}\models q|_{(M,a)}$ and the result follows. 
    
    For the second part, let $J=(g_i)_{i\in\omega}\models q^{\otimes\omega}|_{M}$ and let $a\models p|_{(M,J)}$. Since $p$ is $M$-invariant and $J$ is $M$-indiscernible, also $J$ is $(M,a)$-indiscernible. So again by Lemma \ref{tf_commutes_with_ind_conjugates} we have $[a,a^{g_1g_0^{-1}}]=1$. Also, again by right f-genericity, $g_1g_0^{-1}\models q|_M$, and now since $a\models p|_{(M,g_1g_0^{-1})}$ the result follows.
\end{proof}

Suppose that $q(x)\in S_{G^{00}}(\mathfrak{C})$ is bi-f-generic for $G$ with non-torsion realizations and let $(a,b)\models q^{\otimes 2}|_M$. By Lemma \ref{tf_f-gen_is_2-engel}, we have $[a,a^{b}]=1$ and $[b,b^{a}]=1$. In other words, $[a,b]\in C_G(a)\cap C_G(b)$, so that $\langle a,b\rangle$ is nilpotent of class at most $2$. Thus, by Lemma \ref{gen_n-nil_implies_gen_N_nil} applied with $n=2$ and $N=3$, there is some $\sigma\in S_3$ such that $q^{\otimes 3}(x_1,x_2,x_3)\vdash [x_{\sigma(1)},x_{\sigma(2)},x_{\sigma(3)}]=1$. The following comes by combining refinements of this with the Hall-Witt identity:

\begin{theorem}\label{gen_tf_implies_nil_2}
    Suppose there is an $M$-invariant bi-f-generic type for $G$ concentrated on $G^{00}$ whose realizations are non-torsion. Then there is a definable subgroup $O\leqslant\tilde{G}$, nilpotent of class at most $2$, such that $O\cap G$ has finite index in $G$.
\end{theorem}
\begin{proof}
    Let $q(x)\in S_{G^{00}}(\mathfrak{C})$ be bi-f-generic for $G$ with non-torsion realizations, and let $(a_i)_{i\in\omega}+(b_i)_{i\in\omega}$ be a Morley sequence of $q$ over $M$. By the argument right above this lemma, it follows from Lemma \ref{tf_f-gen_is_2-engel} that $\langle a_1,b_1\rangle$ is nilpotent of class at most $2$. Now, by dp-minimality, one of $(a_i)_{i\in\omega}$ and $(b_i)_{i\in\omega}$ is indiscernible over $[a_1,b_1]$. If $(b_i)_{i\in\omega}$ is, then, since $[a_1,b_1,b_1]=1$, we have $[a_1,b_1,b_2]=1$, so the result follows by Lemma \ref{generic_nil_is_nil}. So we may assume instead that $(a_i)_{i\in\omega}$ is indiscernible over $[a_1,b_1]$. Since $[a_1,b_1,a_1]=1$, thus also $[a_1,b_1,a_2]=1$; since $q$ is bi-f-generic and concentrated on $G^{00}$, and $a_2\models q|_{(M,a_1)}$, we have $a_2^{a_1}\models q|_{(M,a_1)}$, and since $b_1\models q|_{(M,a_1,a_2)}$ it follows that $[a_1,b_1,a_2^{a_1}]=1$.
    
    Now, again by dp-minimality, one of $(a_i)_{i\in\omega}$ and $(b_i)_{i\in\omega}$ is indiscernible over $[b_1,a_1^{b_1^{-1}}]$. If $(b_i)_{i\in\omega}$ is, then, since $[b_1,a_1^{b_1^{-1}},b_1]=1$, also $[b_1,a_1^{b_1^{-1}},b_2]=1$, so that $[b_1,a_1,b_2^{b_1}]=1$. Since $q$ is bi-f-generic for $G$ and concentrated on $G^{00}$, and $b_2\models q|_{(M,a_1,b_1)}$, we have $b_2^{b_1}\models q|_{(M,a_1,b_1)}$, so it follows that $[b_1,a_1,b_2]=1$. But then $[a_1,b_1,b_2]=1$ and so again by Lemma \ref{generic_nil_is_nil} we are done. So we may assume that $(a_i)_{i\in\omega}$ is indiscernible over $[b_1,a_1^{b_1^{-1}}]$. Since $[b_1,a_1^{b_1^{-1}},a_1]=1$, it follows that $[b_1,a_1^{b_1^{-1}},a_0]=1$, and hence that $[b_1,a_1,a_0^{b_1}]=1$.

    Altogether we have shown that, if the desired conclusion does not hold, then $[a,c,b^a]=1$ and $[c,b,a^c]=1$ for all $(a,b,c)\models q^{\otimes 3}|_M$. By the Hall-Witt identity, we have $$[c,b,a^c][a,c,b^a][b,a,c^b]=1,$$ so it follows that $[b,a,c^b]=1$. By one final application of bi-f-genericity of $q$, we have $c^b\models q|_{(M,a,b)}$, so that $[b,a,c]=1$, whence $[a,b,c]=1$ and the result follows from Lemma \ref{generic_nil_is_nil}.
\end{proof}

\begin{theorem}\label{tf_implies_abelian}
    If $G$ is torsion-free, then it is abelian.
\end{theorem}
\begin{proof}
By Fact \ref{dp-min_implies_def_am} and Fact \ref{bi-f-gens_exist}, let $q(x)\in S_{G^{00}}(\mathfrak{C})$ be an $M$-invariant type bi-f-generic for $G$, and let $(a,b)\models q^{\otimes 2}|_M$. As in Theorem \ref{gen_tf_implies_nil_2}, we have that $\langle a,b\rangle$ is nilpotent of class at most $2$, so that $[a^k,b]=[a,b]^k=[a,b^k]$ for all $k\in\omega$. On the other hand, by Fact \ref{dp-min_implies_ab_by_tor}, for some $n>0$, one of $[a^n,b]=1$ and $[a,b^n]=1$ holds, whence $[a,b]^n=1$. Since $G$ is torsion-free, thus $[a,b]=1$, and by Lemma \ref{generic_nil_is_nil} there is a definable abelian subgroup $A\leqslant\tilde{G}$ such that $G\cap A$ has finite-index in $G$. In particular $G^{00}\leqslant G\cap A$ and so $G^{00}$ is abelian.

Now, to show the desired result, it suffices by Fact \ref{tf_central_by_finite} to show that $C_G(g)$ has finite index in $G$ for all $g\in G$. Thus fix arbitrary $g\models G$ with $g\neq 1$. Let $c\models q|_{(M,g)}$. By Lemma \ref{tf_f-gen_is_2-engel} we have $[g,g^c]=1$. On the other hand, $c\models G^{00}$, and $G^{00}$ is a normal abelian subgroup of $G$, so also $[c,c^g]=1$. Thus $\langle c,g\rangle$ is nilpotent of class at most $2$, and so again since $G$ is torsion-free and by Fact \ref{dp-min_implies_ab_by_tor} it follows that $[c,g]=1$. So $c\in C_{\tilde{G}}(g)$. Since $c\models q|_{(M,g)}$ thus $q(y)\vdash y\in C_{\tilde{G}}(g)$; by f-genericity of $q$ it follows from Lemma \ref{f-generic invariant subgroup contains G00} and Fact \ref{bounded index rel definable implies finite index} that $C_G(g)=C_{\tilde{G}}(g)\cap G$ has finite-index in $G$, as needed.
\end{proof}

 \section{Distal case}\label{distal section}
In this section, we will study dp-minimal groups admitting a distal f-generic type; by Corollary \ref{dp-min_implies_gs_or_distal_group_version}, these are precisely those dp-minimal groups with no generically stable f-generic type. Hence the class of such groups includes not only distal dp-minimal groups, ie dp-minimal groups in which every invariant type is distal, but also, for instance, dp-minimal $\mathcal{P}$-valued groups whenever $(\mathcal{P},\leqslant)$ is a total order with no least element.

Many of the results we will prove here were proved for C-minimal valued groups in \cite{simonetta_nilpotent}. In particular, the results of Section 5.3 and Theorem \ref{distal_implies_virtually_nil} were both proved in the C-minimal case there. The main difference here is that we do not start with a valued group structure on the group.

Throughout this section, \uline{assume that $T$ is NIP, that $\tilde{G}$ is a definable group of $T$, and that $G$ is a type-definable dp-minimal subgroup of $\tilde{G}$ that admits a distal f-generic type}. Additionally, by adding constant symbols, for notational convenience \uline{assume that $G$ is type-definable without parameters}. By Fact \ref{bi-f-gens_exist}, and adding constant symbols if necessary, throughout the section \uline{fix a distal type $q(x)\in S_{G^{00}}(\mathfrak{C})$ invariant over any small model and bi-f-generic for $G$}.

\subsection{Conjugacy classes are thin}
Following \cite{gismatullin}, given a subset $X\subseteq G$ with $1\in X=X^{-1}$ and invariant over some small set of parameters, call $X$ `thick' if, for any sequence $(a_i)_{i\in\omega}$ in $G$ indiscernible over the parameters over which $X$ is invariant, we have $a_i^{-1}a_j\in X$ for all $i,j\in\omega$. In this section we will show that, for any $u\in G\setminus\{1\}$, the set $G\setminus (u^G\cup u^{-G})$ is thick.
\begin{lemma}\label{[g,u]_not_equiv_to_u}
    Let $B\subset\mathfrak{C}$ be a parameter set. Suppose $u,g\models G$ are such that $\tp(u/B)\perp^w\tp(g/B)$ and that $\tp(u/B)\perp^w\tp(u^g/B)$, and that there is some $n>0$ with $[g^n,u]=1$. If $[g,u]\equiv_B u$, then $u=1$.
\end{lemma}
\begin{proof}
    Suppose $[g,u]\equiv_B u$. We prove by induction on $k$ that $[g^k,u]\equiv_B u$ for each $k>0$; taking $k=n$ will then give the result. The base case $k=1$ is by hypothesis. For the inductive step, assume that $[g^k,u]\equiv_B u$. Then, since $\tp(u/B)\perp^w\tp(g/B)$, it follows that $[g^k,u]^g\equiv_{B}u^g$. By the case $k=1$, we also have $[g,u]\equiv_B u$. Now since $\tp(u/B)\perp^w\tp(u^g/C)$, the $2$-tuple $([g,u],u^g)$ has the same type over $B$ as $([g,u],[g^k,u]^g)$ does. But we have $u^g\cdot [g,u]=u$, and so the result follows from the identity $[g^k,u]^g[g,u]=[g^{k+1},u]$.
\end{proof}

\begin{lemma}\label{u_perp_ug}
    Let $M$ be a small model, and suppose that $v\in G(M)$ and that $u\in v^{G^{00}}$. If $J+g\models q^{\otimes(\omega+1)}|_{(M,u)}$, then $\tp(u/M,J)\perp^w\tp(g/M,J)$ and $\tp(u/M,J)\perp^w\tp(u^g/M,J)$.
\end{lemma}
\begin{proof}
    That $\tp(u/M,J)\perp^w\tp(g/M,J)$ follows immediately from distality of $q$. Now, let $h\models G^{00}$ be such that $u=v^h$, and let $J'+g'\models q^{\otimes(\omega+1)}|_{(M,u,h)}$. We have $(u,J',g')\equiv_M (u,J,g)$, so it suffices to show that $\tp(u/M,J')\perp^w\tp(u^{g'}/M,J')$. Since $h\models G^{00}$ and $g'\models q|_{(M,u,h,J')}$, it follows by left f-genericity of $q$ that $hg'\models q|_{(M,u,h,J')}$. Since $v\in G(M)$, we thus have in particular that $v^{g'}\equiv_{(M,J')}v^{hg'}=u^{g'}$. On the other hand, again since $v\in G(M)$, and since $\tp(u/M,J')\perp^w\tp(g'/M,J')$, we have $\tp(u/M,J')\perp^w\tp(v^{g'}/M,J')$, so the result follows.
\end{proof}

\begin{lemma}\label{g00_ccs_are_thin}
Let $M$ be a small model, and suppose that $I=(a_i)_{i\in\omega}$ is an $M$-indiscernible sequence of elements of $G$ such that, for some $v\in G(M)$, we have $a_i^{-1}a_j\in v^{G^{00}}$ for all $i<j$. Then $I$ is constant.
\end{lemma}
\begin{proof}
By Ramsey and compactness, it suffices to show the result when $I$ has order type $\mathbb{R}$, so replace $I$ by a sequence $(a_i)_{i\in\mathbb{R}}$ with the same properties. Now, I claim that, for some $n>0$, $q(x)\vdash [x^n,a_i^{-1}a_j]=1$ for all $i,j\in\mathbb{R}$; call this ($\star$). If $q(x)\vdash x^n=1$ for some $n>0$ then the claim is clear. Otherwise the realizations of $q(x)$ are non-torsion, whence by Theorem \ref{gen_tf_implies_nil_2} $G^{00}$ is nilpotent of class at most $2$; in particular then $[b^n,c]=[b,c^n]$ for all $b,c\models G^{00}$ and all $n\in\omega$. Since $a_i^{-1}a_j\models G^{00}$ and $q$ is concentrated on $G^{00}$, ($\star$) then follows from Fact \ref{dp-min_implies_ab_by_tor}.

Let $J+g\models q^{\otimes(\omega+1)}|_{(M,I)}$, and let $u=a_0^{-1}a_1$. By Lemma \ref{u_perp_ug}, we have $\tp(u/M,J)\perp^w\tp(g/M,J)$ and $\tp(u/M,J)\perp^w\tp(u^g/M,J)$; moreover, by \ref{u_perp_ug} applied to the bi-f-generic distal type $q^{-1}$, we also have $\tp(u/M,J)\perp^w\tp(g^{-1}/M,J)$ and $\tp(u/M,J)\perp^w\tp(u^{g^{-1}}/M,J)$. Call these facts ($\ast$).

Note also that, since $I$ is $M$-indiscernible and $q$ is $M$-invariant, $I$ is $(M,J,g)$-indiscernible. Thus $u\equiv_{(M,J,g)}a_i^{-1}a_j$ for all $i<j\in\mathbb{R}$.

Now, $(a_i)_{i\in[0,1]}$ and $(a_j)_{j\in [3,4]}$ are mutually indiscernible over $(M,a_2,a_5,J,g)$. By dp-minimality, one of them is thus indiscernible over $(M,a_2,a_5,J,g,a_0ga_3^{-1})$. First suppose $(a_i)_{i\in[0,1]}$ is indiscernible over that set. Note that $$g^{-1}\cdot a_0^{-1}\cdot (a_0ga_3^{-1})\cdot a_5=a_3^{-1}a_5\equiv_{(M,J)}u.$$ Hence, by indiscernibility, also $g^{-1}a_1^{-1}a_0ga_3^{-1}a_5\equiv_{(M,J)}u$, ie $(a_1^{-1}a_0)^g(a_3^{-1}a_5)\equiv_{(M,J)}u$. But we have $a_1^{-1}a_0\equiv_{(M,J,g)}u^{-1}$, hence $(a_1^{-1}a_0)^g\equiv_{(M,J)} u^{-g}$, and $a_3^{-1}a_5\equiv_{(M,J)}u$. Since $\tp(u^{-g}/M,J)\perp^w\tp(u/M,J)$, it follows that $u^{-g}u\equiv_{(M,J)}u$, ie that $[g,u]\equiv_{(M,J)}u$. By Lemma \ref{[g,u]_not_equiv_to_u}, using ($\ast$) and ($\star$), we must have $u=1$, so that indeed $I$ is constant.

Now suppose $(a_j)_{j\in [3,4]}$ is indiscernible over $(M,a_2,a_5,J,g,a_0ga_3^{-1})$. As above, note that $$a_2^{-1}\cdot (a_0ga_3^{-1})\cdot a_3\cdot g^{-1}=a_2^{-1}a_0\equiv_{(M,J)}u^{-1}.$$ By indiscernibility, thus $a_2^{-1}a_0ga_3^{-1}a_4g^{-1}\equiv_{(M,J)}u^{-1}$, ie $(a_2^{-1}a_0)(a_3^{-1}a_4)^{g^{-1}}\equiv_{(M,J)}u^{-1}$. But we have $a_3^{-1}a_4\equiv_{(M,J,g)}u$, hence $(a_3^{-1}a_4)^{g^{-1}}\equiv_{(M,J)}u^{g^{-1}}$, and $a_2^{-1}a_0\equiv_{(M,J)}u^{-1}$. Since $$\tp(u^{g^{-1}}/M,J)\perp^w\tp(u/M,J),$$ it follows that $u^{-1}u^{g^{-1}}\equiv_{(M,J)}u^{-1}$, so that $u^{-g^{-1}}u\equiv_{(M,J)}u$, ie $[g^{-1},u]\equiv_{(M,J)}u$. By Lemma \ref{[g,u]_not_equiv_to_u} applied to the bi-f-generic distal type $q^{-1}$, using ($\ast$) and ($\star$), we must have $u=1$, so that again $I$ is constant.
\end{proof}

\begin{corollary}\label{ccs_are_thin}
Suppose $I=(b_i)_{i\in\omega}$ is an indiscernible sequence of elements of $G$ such that $b_i^{-1}b_j\sim_G b_0^{-1}b_1$ for all $i<j$. Then $I$ is constant.
\end{corollary}
\begin{proof}
By Ramsey and compactness, we may find a model $M$ such that (i) $G(M)$ contains a representative of every coset of $G^{00}$, (ii) $I$ is $M$-indiscernible, and (iii) for some $v\in G(M)$, we have $b_i^{-1}b_j\in v^G$ for all $i<j$. By (i), there is some $\lambda\in G(M)$ such that $b_i^{-1}b_j\in v^{G^{00}\lambda^{-1}}$ for all $i<j$. Let $a_i=b_i^\lambda$. Then $(a_i)_{i\in\omega}$ is $M$-indiscernible, and $a_i^{-1}a_j\in v^{G^{00}}$ for all $i<j$; by Lemma \ref{g00_ccs_are_thin}, $(a_i)_{i\in\omega}$ must be constant, so $I$ is too.
\end{proof}

By considering a Morley sequence we hence obtain the following (note that the relation $x\sim_G y$ is type-definable):

\begin{corollary}\label{ccs_are_thin_type_version}Let $M\prec\mathfrak{C}$, and suppose $p(x)\in S_G(\mathfrak{C})$ is $M$-invariant with $p^{\otimes 2}(x,y)\vdash p(x^{-1}y)$ and $p^{\otimes 2}(x,y)\vdash x\sim_G y$. Then $p(x)\vdash x=1$.
\end{corollary}

\subsection{Conjugacy classes determine types}
In this section we will show that any type concentrated on $G$ but not concentrated on boundedly many conjugacy classes of $G$ is uniquely determined by the $G^{00}$-conjugacy class of any of its realizations.
\begin{lemma}\label{cc_determines_types_seq_form}
    Suppose $B\subset\mathfrak{C}$ is a small set of parameters and that $a\models G$ begins a $B$-indiscernible sequence of elements that all lie in different conjugacy classes of $G$. Then $a^g\equiv_B a$ for all $g\models G^{00}$.
\end{lemma}
\begin{proof}
    By Ramsey and compactness, we may assume without loss of generality that $B=M$ is a model. I first claim that, for any $J+g\models q^{\otimes(\omega+1)}|_{(M,a)}$, we have $a^g\equiv_{(M,J)} a$; it suffices for this to show that $a^g\equiv_{(M,J)} a$ for some $J+g\models q^{\otimes(\omega+1)}|_{(M,a)}$. To see this, let $I=(a_i)_{i\in\omega}$ be an $M$-indiscernible sequence with $a_0=a$ and $a_i\not\sim_G a_j$ for all $i\neq j$. Let $J\models q^{\otimes\omega}|_{(M,I)}$ and let $K=(g_i)_{i\in\omega}\models q^{\otimes\omega}|_{(M,I,J)}$. Since $I$ is $M$-indiscernible and $q$ is $M$-invariant, $I$ and $K$ are mutually indiscernible over $(M,J)$, and thus by dp-minimality one of them is indiscernible over $(M,J,a^{g_1})$. Since $a_i\sim_G a^{g_1}=a_0^{g_1}$ if and only if $i=0$, $K$ must therefore be indiscernible over $(M,J,a^{g_1})$, and hence there is some $a_*\equiv_{(M,J)} a$ such that $a_*^{g_0}=a^{g_1}$. But now $a^{g_1g_0^{-1}}=a_*$, so that $a^{g_1g_0^{-1}}\equiv_{(M,J)} a$. Since $q$ is right f-generic and concentrated on $G^{00}$, and $g_1\models q|_{(M,a,J,g_0)}$, we have $g_1g_0^{-1}\models q|_{(M,a,J)}$, so the claim follows.

    Now, let $J+g\models q^{\otimes(\omega+1)}|_{(M,a)}$. By the claim, $a^g\equiv_{(M,J)} a$, and by the claim applied to the bi-f-generic distal type $q^{-1}$ also $a^{g^{-1}}\equiv_{(M,J)}a$. On the other hand, by distality of $q$, we have $\tp(a/M,J)\perp^w\tp(g/M,J)$, so it follows that $a_*^g\equiv_{(M,J)}a\equiv_{(M,J)}a_*^{g^{-1}}$ for all $a_*\equiv_{(M,J)}a$. In other words, if $r(x)=\tp(a/M,J)$, then $g$ lies in the normalizer $N_G(r):=\{h\in G:r(x)^h=r(x)\}$. $N_G(r)$ is an $(M,J)$-invariant subgroup of $G$; since $\tp(g/M,J)=q|_{(M,J)}$ and $q$ is f-generic, it follows by Fact \ref{f-generic invariant subgroup contains G00} that $G^{00}\leqslant N_G(r)$, so that $a^h\equiv_{(M,J)}a$ for all $h\models G^{00}$. The result follows `a fortiori'.
\end{proof} The following corollary is standard.
\begin{corollary}\label{cc_determines_types}
    If $r(x)\in S_G(B)$ for some $B\subset\mathfrak{C}$ and $r$ is not concentrated on boundedly many conjugacy classes of $G$, then $r(x)^g=r(x)$ for all $g\models G^{00}$. If $B=M$ is a model, then the same holds provided $r$ is not concentrated on a single conjugacy class.
\end{corollary}
\begin{proof}
    For the first part, from a long sequence of realizations of $r$ lying in different conjugacy classes of $G$ we may use Erd\H{o}s-Rado to obtain a $B$-indiscernible sequence of realizations of realizations of $r$ lying in different conjugacy classes of $G$, and then Lemma \ref{cc_determines_types_seq_form} applies. The second part follows from the fact that Lascar strong types coincide with types over models. Indeed the relation $x\sim_G y$ is a type-definable equivalence relation on $G$, and by the first part it has boundedly many classes in $r$. So by a standard argument (eg Proposition 5.11 in \cite{simon_book}) there is an $M$-invariant bounded equivalence relation whose restriction to $r$ coincides with $\sim_G$. But since $M$ is a model $r$ implies a complete Lascar strong type over $M$, hence concentrates on a single equivalence class of $\sim_G$.
\end{proof}

\begin{corollary}\label{cc_determines_morley_seq}
Suppose $p(x)\in S_G(\mathfrak{C})$ is an $M$-invariant type not concentrated on a conjugacy class of $G$. If $(a_i)_{i\in\omega}\models p^{\otimes\omega}|_M$, then $(a_i^{g_i})_{i\in\omega}\models p^{\otimes\omega}|_M$ for all $g_i\models G^{00}$.
\end{corollary}
\begin{proof}
If $(b_i)_{i\in \omega}$ is a Morley sequence of $p$ over $M$, then each $b_i$ starts an $(M,b_{\neq i})$-indiscernible sequence of non-conjugates, so by Lemma \ref{cc_determines_types_seq_form} $b_i^g\equiv_{(M,b_{\neq i})}b_i$ for all $g\models G^{00}$. Thus by induction on $n\in\omega$ we have $(a_i^{g_i})_{i<n}+(a_i)_{i\geqslant n}\models p^{\otimes\omega}|_M$ for all $g_i\models G^{00}$, and the result follows. \end{proof}

\subsection{Commutators}
In this section we will show that, for any $a\models G$, the set of commutators $[a,G^{00}]$ is a group; this was proved for C-minimal valued groups in \cite{simonetta_nilpotent}.
\begin{lemma}\label{[a,q]_not_concentrated_on_cc}
Let $a\models G$ such that $q(x)\vdash [a,x]\neq 1$. Then $[a,q]$, ie the pushforward of $q$ under the map $x\mapsto [a,x]$, is not concentrated on a single conjugacy class of $G$.
\end{lemma}
\begin{proof}
Suppose otherwise for contradiction that $[a,q](x)\vdash x\sim_G u$ for some $u\models G$. Let $M$ be any small model containing $a,u$, and let $(g_i)_{i\in\omega}\models q^{\otimes\omega}|_M$. Since $q$ is right f-generic and concentrated on $G^{00}$, we have $g_jg_i^{-1}\models q|_M$ for all $i<j$ and hence $[a,g_jg_i^{-1}]\sim_G u$ for all $i<j$. For $i\in\omega$, let $b_i=a^{g_i}$, and consider the $M$-indiscernible sequence $I=(b_i)_{i\in\omega}$. For $i<j$, we have $$b_i^{-1}b_j=a^{-g_i}a^{g_j}=[a,g_jg_i^{-1}]^{g_i}\sim_G u.$$ By Lemma \ref{ccs_are_thin}, $I$ is constant, so that $a^{g_0}=a^{g_1}$ and hence $g_1g_0^{-1}\in C_G(a)$. Since $g_1g_0^{-1}\models q|_M$, this contradicts that $q(x)\vdash [a,x]\neq 1$.
\end{proof}

By Corollary \ref{cc_determines_types}, and since $[a,q]$ is just $\tp(1/\mathfrak{C})$ if $q(x)\vdash [a,x]=1$, we get the following:

\begin{corollary}\label{[a,q|_N]_is_normal}
For any $M\prec\mathfrak{C}$ and any $a\in G(M)$, we have $[a,q|_M]=[a,q|_M]^g$ for all $g\models G^{00}$.
\end{corollary}

\begin{lemma}\label{[a,q]_is_[a,G]-inv}
Let $a\models G$. Then $[a,g][a,q]=[a,q]$ for all $g\models G^{00}$, and $[a,G^{00}][a,q|_M]\subseteq [a,G^{00}]$ for any small model $M$ with $a\in G(M)$.
\end{lemma}
\begin{proof}
First note the following general fact. Let $r\in S_G(M)$ be a complete type over some model $M\prec\mathfrak{C}$ such that $a\in G(M)$ and such that $[a,r]=[a,r]^g$ for all $g\models G^{00}$. Then, from the identity $[a,g][a,c]^g=[a,cg]$, it follows that $[a,g][a,r]=[a,rg]$ for all $g\models G^{00}$.

The desired result now follows; indeed by Corollary \ref{[a,q|_N]_is_normal} and by the remark above we have $[a,g][a,q|_M(x)]=[a,q|_M(x)g]$ for all $M\prec\mathfrak{C}$ with $a\in G(M)$ and all $g\models G^{00}$; this shows already that $[a,G^{00}][a,q|_M]\subseteq [a,G^{00}]$. Furthermore, since $q$ is right f-generic, we have $qg=q$ for all $g\in G^{00}(\mathfrak{C})$, so that $q|_M(x)g=q|_M(x)$ whenever $G(M)\ni g$, and hence indeed $[a,g][a,q]=[a,q]$.
\end{proof}

\begin{lemma}\label{[a,G]_is_group} For any $a\models G$, the set of commutators $[a,G^{00}]$ is a type-definable subgroup of $G^{00}$, and it is `strongly connected', ie has no proper type-definable subgroups of bounded index.
\end{lemma}
\begin{proof} Certainly $[a,G^{00}]$ is type-definable, and it is a subset of $G^{00}$ since $G^{00}$ is  a normal subgroup of $G$. Now, let $M$ be any small model containing $a$. By Corollary \ref{[a,q|_N]_is_normal}, we have $[a,q|_M]^g=[a,q|_M]$ for all $g\models G^{00}$; likewise, by Corollary \ref{[a,q|_N]_is_normal} applied to the bi-f-generic distal type $q^{-1}$, we also have $[a,q|_M^{-1}]^g=[a,q|_M^{-1}]$ for all $g\models G^{00}$. Thus, for any $b\models q|_M$, since $q$ is concentrated on $G^{00}$ we have
\begin{align*}
[a,b^{-1}]&=[a,b]^{-b^{-1}}\in [a,q|_M]^{-b^{-1}}=[a,q|_M]^{-1} \text{ and} \\
[a,b]^{-1}&=[a,b^{-1}]^{b^{-1}}\in [a,q|_M^{-1}]^{b^{-1}}=[a,q|_M^{-1}],
\end{align*} so that $[a,q|_M^{-1}]=[a,q|_M]^{-1}$, and it follows that $[a,q|_M][a,q|_M^{-1}]$ is closed under taking inverses. Now, by Fact \ref{f-generic facts}(\ref{f-gen_iff_g00_inv}), we have $G^{00}\subseteq q|_M^{-1}(x)q|_M(x)$. Thus, for any $g\models G^{00}$, there are $b,c\models q|_M$ with $g=b^{-1}c$, and then we have $$[a,g]=[a,c][a,b^{-1}]^c\in [a,q|_M][a,q|_M^{-1}]^c=[a,q|_M][a,q|_M^{-1}],$$ so that $[a,G^{00}]\subseteq [a,q|_M][a,q|_M^{-1}]$. Finally, by Lemma \ref{[a,q]_is_[a,G]-inv} applied to $q$ and to $q^{-1}$, we have $$[a,G^{00}][a,q|_M][a,q|_M^{-1}]\subseteq [a,G^{00}][a,q|_M^{-1}]\subseteq [a,G^{00}].$$ Thus $[a,q|_M][a,q|_M^{-1}]=[a,G^{00}]$, and it is symmetric and closed under multiplication; in other words, it is a group.

To see that $[a,G^{00}]$ is strongly connected, note that, by Lemma \ref{[a,q]_is_[a,G]-inv}, $[a,q]$ is invariant under left translation by all elements of $[a,G^{00}(\mathfrak{C})]$. If $H\leqslant [a,G^{00}]$ is type-definable and of bounded index in $[a,G^{00}]$, then, since $[a,q]$ is concentrated on $[a,G^{00}]$, it must be concentrated on some right coset of $H$, and then left-invariance of $[a,q]$ forces $H=[a,G^{00}]$.
\end{proof}

\begin{corollary}\label{[a,G]s_are_chain}
For any $a,b\models G$, one of $[a,G^{00}]\supseteq[b,G^{00}]$ and $[a,G^{00}]\subseteq [b,G^{00}]$ holds.
\end{corollary}
\begin{proof} By Fact \ref{inv_subgps_are_commensurable}, $[a,G^{00}]\cap [b,G^{00}]$ has bounded index in one of $[a,G^{00}]$ and $[b,G^{00}]$. Since the latter two groups are strongly connected, $[a,G^{00}]\cap [b,G^{00}]$ must be equal to one of them, and the result follows.
\end{proof}

\begin{corollary}
    For any type-definable set $X\subseteq G$, the set of commutators $[X,G^{00}]$ is a type-definable strongly connected subgroup of $G^{00}$.
\end{corollary}
\begin{proof}
    $[X,G^{00}]$ is equal to $\bigcup_{a\models X}[a,G^{00}]$, and by Lemma \ref{[a,G]_is_group} and Corollary \ref{[a,G]s_are_chain} this is the union of a chain of subgroups of $G^{00}$. So indeed $[X,G^{00}]$ is a group, and it is certainly type-definable. To see that it is strongly connected, suppose that $H\leqslant [X,G^{00}]$ is a bounded-index type-definable subgroup. Then, for every $a\models X$, $H\cap [a,G^{00}]$ has bounded index in $[a,G^{00}]$, hence must be equal to $[a,G^{00}]$. So $H\geqslant [a,G^{00}]$ for all $a\models X$ and hence $H=[X,G^{00}]$, as desired.
\end{proof}

The following now has an identical proof as Proposition 2.4.1 of \cite{simonetta_nilpotent}:

\begin{lemma}\label{[G00,G00]_has_fin_exp}
    $[G^{00},G^{00}]$ has finite exponent.
\end{lemma}
\begin{proof}
    Fix $a,b\models G^{00}$; by Lemma \ref{[a,G]_is_group}, $[a,b,G^{00}]$ is now a normal subgroup of $G^{00}$. Working in the quotient $G^{00}/[a,b,G^{00}]$, it follows from the identities $[a^{k+1},b]=[a,b][a,b,a^k][a^k,b]$ and $[a,b^{k+1}]=[a,b^k][a,b][a,b,b^k]$ that $[a^k,b]$ and $[a,b^k]$ are each equivalent to $[a,b]^k$ modulo $[a,b,G^{00}]$ for each $k\in\omega$. By Fact \ref{dp-min_implies_ab_by_tor}, there is some $n>0$ such that one of $[a^n,b]$ and $[a,b^n]$ is trivial, and it follows that $[a,b]^n\in [a,b,G^{00}]$. 

    Thus there is some $g\models G^{00}$ with $[a,b]^n=[a,b,g]$, so that $[a,b]^{n+1}=[a,b]^g\sim [a,b]$. By Lemma \ref{conjugate_with_power_in_ab_by_torsion} and Fact \ref{dp-min_implies_ab_by_tor}, it follows that $[a,b]$ must be a torsion element. Since $a,b\models G^{00}$ were arbitrary, thus every element of $[G^{00},G^{00}]$ has finite order, and so by compactness $[G^{00},G^{00}]$ has finite exponent.
\end{proof}

\begin{corollary}\label{distal_implies_virtually_nil_2_or_torsion}
    If $G$ does not have finite exponent, then it is nilpotent-by-finite of class $\leqslant 2$.
\end{corollary}
\begin{proof}
    If the realizations of $q(x)$ are non-torsion, then by Theorem \ref{gen_tf_implies_nil_2} we are done, so we may assume that $q(x)\vdash x^k=1$ for some $k>0$; we aim to show that $G$ has finite exponent. By Lemma \ref{[G00,G00]_has_fin_exp}, let $m>0$ be the exponent of $[G^{00},G^{00}]$. By Fact \ref{f-generic facts}(\ref{f-gen_iff_g00_inv}), every element of $G^{00}$ is a difference of two realizations of $q$, and it follows that $G^{00}$ has exponent at most $mk$: if $a,b\models G^{00}$ have order at most $k$, then $(ab)^k$ is equivalent to $a^kb^k=1$ modulo $[G^{00},G^{00}]$, so that $(ab)^k$ lies in $[G^{00},G^{00}]$, and hence $(ab)^{km}=1$. Thus indeed $G^{00}$ has exponent at most $n:=mk$.

    Now, since $G^{00}$ is a group, we have $G^{00}(x)\wedge G^{00}(y)\vdash (xy^{-1})^n=1$, and so by compactness there is some $\varnothing$-definable set $D\supseteq G^{00}$ such that $(ab^{-1})^n=1$ for all $a,b\models D$; in other words $c^n=1$ for all $c\in DD^{-1}$.

    Since $G$ does not have finite exponent, the partial type $G(y)\wedge\{y^n\neq 1:n>0\}$ is finitely satisfiable, and so by compactness there is a non-torsion element $b\models G$. Since every element of $DD^{-1}$ is torsion, we have $b^k\notin DD^{-1}$ for all $k\neq 0$, and hence $b^iD\cap b^jD=\varnothing$ for all $i\neq j\in\omega$. Letting $I=(a_i)_{i\in\omega}$ be an indiscernible sequence realizing the EM-type of $(b^i)_{i\in\omega}$ over $\varnothing$, we then have $a_iD\cap a_jD=\varnothing$ for all $i\neq j$; since $D\supseteq G^{00}$, the $a_i$ thus lie in pairwise distinct cosets of $G^{00}$. Since $I$ is indiscernible, this contradicts that $G^{00}$ has bounded index in $G$, and the result follows.
\end{proof}

% \begin{remark}
%     Note that this fails in infinite dihedral group if we weren't working under assumption x.xx.
% \end{remark}

\subsection{Valued group structure}
Suppose that $(H,\cdot)$ is an arbitrary group and that $N\trianglelefteqslant H$ is a normal subgroup such that the collection $\mathcal{P}=\{[a,N]:a\in H\}$ is a family of subgroups of $H$ totally ordered by the reverse inclusion relation, which we will denote $\preccurlyeq$. Note that $(\mathcal{P},\preccurlyeq)$ has a maximal element $\{1\}=[1,N]$. Now, define a map $v:H\to\mathcal{P}$ by $v(a)=[a,N]$. By definition, $v(a)\preccurlyeq v(b)$ if and only if $[a,N]\supseteq[b,N]$. Note that $v$ is constant on cosets of $C_H(N)$, since $[C_H(N)a,g]=\{[a,g]\}$ for all $a\in H$ and $g\in N$. So $v$ descends to a well-defined map $H/C_H(N)\to \mathcal{P}$, which we will also denote $v$.

\begin{lemma}\label{G/FC_is_valued}
In the above circumstances, $(H/C_H(N),v)$ is a $\mathcal{P}$-valued group.
\end{lemma}
\begin{proof}
For property (1), note that, for $a\in H$, we have $v(a)=\max\mathcal{P}$ if and only if $[a,N]=\{1\}$, ie if and only if $a\in C_H(N)$. For property (2), suppose that $v(a)\prec v(b)$ and that $g\in H$ is arbitrary. We have $[a,N]\supsetneq [b,N]$, and hence $[a,N]^g\supsetneq [b,N]^g$. But $N$ is a normal subgroup of $H$, so $[a,N]^g=[a^g,N^g]=[a^g,N]$, and likewise $[b,N]^g=[b^g,N]$, so that $v(a^g)\prec v(b^g)$ as desired. For property (3), fix $a\in H$. We have $[a^{-1},N]=[N,a]^{a^{-1}}=[N,a]$. But $[a,N]$ is by hypothesis a subgroup of $H$, hence in particular is closed under inverses, and so $[a^{-1},N]=[a,N]$, whence $v(a)=v(a^{-1})$.

So we need only show property (4), ie that $v(ab)\succcurlyeq\min\{v(a),v(b)\}$ for all $a,b\in H$. First suppose that $v(b)=\min\{v(a),v(b)\}$, ie that $[a,N]\subseteq [b,N]$; we wish to show that $[ab,N]\subseteq [b,N]$. By the identity $[ab,g]=[a,g]^b[b,g]$, we have $[ab,N]\subseteq [a,N]^b[b,N]$. But $[a,N]^b\subseteq [b,N]^b=[b,N]$, so we have $[ab,N]\subseteq [b,N][b,N]\subseteq [b,N]$, where the second inclusion uses that $[b,N]$ is closed under multiplication. In other words $v(ab)\succcurlyeq v(b)$, as desired.

Now suppose that $v(a)=\min\{v(a),v(b)\}$. Using that $v(c)=v(c^{-1})$ for all $c\in H$, we have $v(ab)=v(b^{-1}a^{-1})$ and $v(b^{-1})=v(b)\succcurlyeq v(a)=v(a^{-1})$. Thus $v(a^{-1})=\min\{v(b^{-1}),v(a^{-1})\}$ and so, by the case in the previous paragraph, we have $v(b^{-1}a^{-1})\succcurlyeq v(a^{-1})$. But this means $v(ab)\succcurlyeq v(a)$, as desired.
\end{proof}

Now, let us return to the group $G$. The subgroup $N=G^{00}$ is normal in $G$, and by Lemma \ref{[a,G]_is_group} and Corollary \ref{[a,G]s_are_chain} the set $\mathcal{P}=\{[a,G^{00}]:a\models G\}$ is a chain of subgroups of $G$. So, by Lemma \ref{G/FC_is_valued}, the map $v:a\mapsto [a,G^{00}]$ induces a $\mathcal{P}$-valued group structure on $G/C_G(G^{00})$. Note that, by Fact \ref{g00_exists} and Fact \ref{bounded index rel definable implies finite index}, the group $C_G(G^{00})$ has a purely group-theoretic description; indeed recall that the FC-center of $G$, which we will denote $\mathrm{FC}(G)$, is the set of all $g\in G$ with $[G:C_G(g)]<\infty$. This is a characteristic subgroup of $G$, and by Fact \ref{g00_exists} and Fact \ref{bounded index rel definable implies finite index}, using that $C_G(g)=C_{\tilde{G}}(g)\cap G$ is relatively definable, we have $g\in\mathrm{FC}(G)$ if and only if $G^{00}\subseteq C_G(g)$. So $C_G(G^{00})=\mathrm{FC}(G)$ and thus the map $v$ makes $G/\mathrm{FC}(G)$ into a valued group.

We have $\mathrm{FC}(G)\cap G^{00}=C_G(G^{00})\cap G^{00}=Z(G^{00})$. Thus $G^{00}/Z(G^{00})$ inherits the valued group structure from $G/\mathrm{FC}(G)$.

\begin{theorem}\label{g00/Z_is_plain}
$(G^{00}/Z(G^{00}),v)$ is a plain valued group.
\end{theorem}
\begin{proof}
We need only show that $v$ is a plain valuation. Hence fix $a,b\in G^{00}\setminus Z(G^{00})$; we wish to show that $v([a,b])\succ\max\{v(a),v(b)\}$. Since $v([a,b])=v([b,a])$, we may assume without loss that $v(a)=\max\{v(a),v(b)\}$, so that $[a,G^{00}]\subseteq [b,G^{00}]$. We wish to show that $v([a,b])\succ v(a)$, ie that $[a,b,G^{00}]\not\supseteq [a,G^{00}]$. To see this, let $M$ be any small model and let $g\models q|_{(M,a,b)}$; it suffices to show that $[a,g]\notin [a,b,G^{00}]$. Suppose otherwise for contradiction that $[a,g]=[a,b,h]$ for some $h\models G^{00}$. Now, $[a,g]\models [a,q]|_{(M,a,b)}$, and $b\models G^{00}$, so by Lemma \ref{[a,q]_is_[a,G]-inv} we have $[a,b][a,g]\models [a,q]|_{(M,a,b)}$. But $[a,b][a,g]=[a,b][a,b,h]=[a,b]^h$. Thus $[a,q]$ is concentrated on the conjugacy class $[a,b]^G$, and by Lemma \ref{[a,q]_not_concentrated_on_cc} we must have $q(x)\vdash x\in C_G(a)$. But then, by f-genericity of $q$ and Fact \ref{f-generic invariant subgroup contains G00} and Fact \ref{bounded index rel definable implies finite index}, $C_G(a)=C_{\widetilde{G}}(a)\cap G$ has finite index in $G$, contradicting that $a\notin\mathrm{FC}(G)$.
\end{proof}

% \begin{remark}
%     Note that, in the example of a c-minimal group that is not abelian-by-finite, this recovers the valuation on $G/Z(G)$.
% \end{remark}

Now by Fact \ref{plain_implies_loc_nil} we obtain the following:

\begin{corollary}\label{G00_is_loc_nil} $G^{00}$ is locally nilpotent.
\end{corollary}
\begin{proof}
    If $G$ has finite exponent, then by Fact \ref{plain_implies_loc_nil} and Theorem \ref{g00/Z_is_plain} $G^{00}/Z(G^{00})$ is locally nilpotent, whence $G^{00}$ is too. So we may assume that $G$ does not have finite exponent, and then it follows from Corollary \ref{distal_implies_virtually_nil_2_or_torsion} that $G^{00}$ is nilpotent of class at most $2$, hence in particular locally nilpotent.
\end{proof}

\subsection{Local nilpotence to nilpotence}
From Lemma \ref{gen_n-nil_implies_gen_N_nil}, since $G$ has dp-rank 1, we have the following:

\begin{lemma}\label{gen_2-nil_implies_gen_nil}
    Let $M\prec\mathfrak{C}$ be any model, and suppose $p(x)\in S_G(\mathfrak{C})$ is an $M$-invariant type with the property that, for some (any) $(a,b)\models p^{\otimes 2}|_M$, the group $\langle a,b\rangle$ is nilpotent of class at most $(n-1)$. Then there is some $\sigma\in S_n$ such that $p^{\otimes n}(x_1,\dots,x_n)\vdash [x_{\sigma(1)},\dots,x_{\sigma(n)}]=1$.
\end{lemma}

Now we can obtain the main result.

\begin{theorem}\label{distal_implies_virtually_nil}
    $G$ is nilpotent-by-finite.
\end{theorem}
\begin{proof}
    By Corollary \ref{G00_is_loc_nil}, $G^{00}$ is locally nilpotent, so any global type $p(x)$ concentrated on $G^{00}$ and invariant over some small set satisfies the hypothesis of Lemma \ref{gen_2-nil_implies_gen_nil} for some $n>0$. In particular, there is some $n>0$ and some $\sigma\in S_n$ such that $q^{\otimes n}(x_1,\dots,x_n)\vdash [x_{\sigma(1)},\dots,x_{\sigma(n)}]=1$. Fix any sequence $(a_1,\dots,a_n)\models q^{\otimes n}|_M$. Since $q(x)$ is concentrated on $G^{00}$, we have $q^{\otimes 2}(x,y)\vdash q(x^{-1}y)$ by Fact \ref{f-generic facts}(\ref{f-gen_iff_g00_inv}), and so by Corollary \ref{ccs_are_thin_type_version} $q$ is not concentrated on any conjugacy class of $G$. By Corollary \ref{cc_determines_morley_seq}, we thus have $(a_1^{g_1},\dots,a_n^{g_n})\models q^{\otimes n}|_M$ for all $g_i\models G^{00}$, and hence in particular $[a_{\sigma(1)}^{G^{00}},\dots,a_{\sigma(n)}^{G^{00}}]=\{1\}$. For each $k\in\{0\}\cup[n]$ and $l\in[n]$, define $A_{kl}=a_l^{G^{00}}$ if $l\leqslant n-k$ and $A_{kl}=G^{00}$ otherwise. We prove by induction on $k$ that $[A_{k\sigma(1)},\dots,A_{k\sigma(n)}]=\{1\}$; we have already shown the base case $k=0$.

    Now, for the inductive step, suppose we have shown the result for $k$, and let $i=\sigma^{-1}(n-k)$. Note that each $A_{kl}$ is setwise invariant under conjugation by elements of $G^{00}$. Hence, by Lemma \ref{commutator_of_normal_subsets}, the set of $g\models G^{00}$ such that $$[A_{k\sigma(1)},\dots,A_{k\sigma(i-1)},g,A_{k\sigma(i+1)},\dots,A_{kn}]=\{1\}$$ is a subgroup of $G^{00}$; call it $U$. By the inductive hypothesis, $U\supseteq A_{k\sigma(i)}=a_{n-k}^{G^{00}}$, and hence in particular $a_{n-k}\in U$. On the other hand, for $j\neq i$, $A_{k\sigma(j)}$ is either equal to $G^{00}$ or of form $a_l^{G^{00}}$ for some $l<n-k$, and it follows that $U$ is $(M,a_{<n-k})$-invariant. But $a_{n-k}\models q|_{(M,a_{<n-k})}$ and $a_{n-k}\in U$, so all the realizations of $q|_{(M,a_{<n-k})}$ lie in $U$, and by Fact \ref{f-generic invariant subgroup contains G00} it follows that $G^{00}\subseteq U$ and hence $G^{00}=U$. Since $A_{(k+1),\sigma(i)}=G^{00}$ and $A_{(k+1),\sigma(j)}=A_{k\sigma(j)}$ for $j\neq i$, the inductive step follows.

    Taking $k=n$, we get $[G^{00},_{\ n-1}G^{00}]=\{1\}$, and the result follows. (For example, we now have in particular that $q^{\otimes n}(x_1,\dots,x_n)\vdash [x_1,\dots,x_n]=1$, and so the result follows by Lemma \ref{generic_nil_is_nil}.)
\end{proof}

\newpage

%\bibliographystyle{plain}
%\bibliography{references}

\begin{thebibliography}{9}
\bibitem{de_aldama} Ricardo de Aldama. Definable nilpotent and soluble envelopes in groups without the
independence property. \textit{Mathematical Logic Quarterly}, Vol. 59, No. 3 (2013). pp. 201-205.

\bibitem{baldwin_saxl} John Baldwin and Jan Saxl. Logical stability in group theory. \textit{Journal of the Australian Mathematical Society (Series A)}, Vol. 21 No. 3 (1976). pp. 267-276.

\bibitem{bryant_hartley}RM Bryant and B Hartley. Periodic locally soluble groups with the minimal condition on centralizers. \textit{Journal of Algebra}, Vol. 61, Iss. 2 (1979). pp. 328-334.

\bibitem{chernikov_pillay_simon}Artem Chernikov, Anand Pillay, Pierre Simon. External definability and groups in NIP theories. \textit{Journal of the London Mathematical Society}, Vol. 90 Iss. 1 (2014). pp 213-240.

\bibitem{chernikov_simon}Artem Chernikov, Pierre Simon. Definably amenable NIP groups. \textit{Journal of the AMS}, Vol. 31 Iss. 3 (2018). pp. 609-641.

\bibitem{profinite}Tim Clausen. Dp-minimal profinite groups and valuations on the integers. Preprint (2020), arXiv:2008.08797. 

\bibitem{delbee_halevi}Christian d'Elbée and Yatir Halevi. Dp-Minimal integral domains. \textit{Israel Journal of Mathematical}, Vol. 246 (2021).

\bibitem{delbee_halevi_johnson}Christian d'Elbée, Yatir Halevi, Will Johnson. The classification of dp-minimal integral domains. \textit{Model Theory}, Vol. 4 No. 2 (2025).

\bibitem{dobrowolski_wagner}Jan Dobrowolski and Frank Wagner. On $\omega$-categorical groups and rings of finite burden. \textit{Israel Journal of Mathematics}, Vol. 236 (2020). pp. 801-839.

\bibitem{dobrowolski_goodrick}Jan Dobrowolski and John Goodrick. Some remarks on inp-minimal and finite burden groups. \textit{Archive for Mathematical Logic}, Vol. 58 (2019). pp. 267-274.

\bibitem{dolich_goodrick_lippel}Alfred Dolich, John Goodrick, David Lippel. Dp-minimality: basic facts and examples. \textit{Notre Dame Journal of Formal Logic}, Vol. 52 No. 3 (2011). pp. 267-288.

\bibitem{rosy}Clifton Ealy, Krzysztof Krupiński, Anand Pillay. Superrosy dependent groups having finitely satisfiable generics. \textit{Annals of Pure and Applied Logic}, Vol. 151 (2008). pp. 1-21.

\bibitem{gismatullin}Jakub Gismatullin. Model-theoretic connected components of groups. \textit{Israel Journal of Mathematics}, Vol. 184 (2011). pp. 251-274.

\bibitem{guingona}Vince Guingona. On VC-minimal fields and dp-smallness. \textit{Archive for Mathematical Logic}, Vol. 53 (2014).

\bibitem{houcine}Abderezak Ould Houcine. A remark on the definability of the Fitting subgroup and the soluble radical. \textit{Mathematical Logic Quarterly}, Vol. 59, No. 1-2 (2013). pp. 62-65.

\bibitem{hrushovski_peterzil_pillay}Ehud Hrushovski, Ya'acov Peterzil, Anand Pillay. Groups, measures, and the NIP. \textit{Journal of the AMS}, Vol. 21 No. 2 (2008). pp. 563-596.

\bibitem{hrushovski_pillay}Ehud Hrushovski, Anand Pillay. On NIP and invariant measures. \textit{Journal of the EMS}, Vol. 13 No. 4 (2011). pp. 1005-1061.

\bibitem{hrushovski_pillay_simon}Ehud Hrushovski, Anand Pillay, Pierre Simon. A note on generically stable measures and fsg groups. \textit{Notre Dame Journal of Formal Logic}, Vol. 53 No. 4 (2012). pp. 599-605.

\bibitem{jahnke_simon_walsberg}Franziska Jahnke, Pierre Simon, Erik Walsberg. Dp-minimal valued fields. \textit{The Journal of Symbolic Logic}, Vol. 82, No. 1 (2017).

\bibitem{johnson_topology}Will Johnson. The canonical topology on dp-minimal fields. \textit{Journal of Mathematical Logic}, Vol. 18, No. 2 (2018).

\bibitem{johnson}Will Johnson. The classification of dp-minimal and dp-small fields. \textit{Journal of the EMS}, Vol. 25 No. 2 (2022). pp. 467–513.

\bibitem{kaplan_levi_simon}Itay Kaplan, Elad Levi, Pierre Simon. Some remarks on dp-minimal groups, in: Groups, Modules, and Model Theory, Surveys and Recent Developments. \textit{Springer} (2017).

\bibitem{krupinski_portillo}Krzysztof Krupiński and Adrián Portillo. On stable quotients. \textit{Notre Dame Journal of Formal Logic}, Vol. 63 No. 3 (2022). pp. 373-394.

\bibitem{onshuus_usvyatsov}Alf Onshuus and Alexander Usvyatsov. On dp-minimality, strong dependence, and weight. \textit{Journal of Symbols Logic}, Vol. 71 No. 1 (2006). pp. 1-21.

\bibitem{pillay_tanovic}Anand Pillay and Predrag Tanović. Generic stability, regularity, and quasi-minimality. Preprint (2009). arXiv:0912.1115.

\bibitem{reineke}Joachim Reineke. Minimale Gruppen. \textit{Z. Math. Logik Grundlagen Math}, Vol. 21 No. 4 (1975). pp. 357–359.

\bibitem{robinson}Derek Robinson. A Course in the Theory of Groups. \textit{Springer} (1996).

 \bibitem{shelah_gs}Saharon Shelah. Classification theory for theories with NIP - a modest beginning, \textit{Scientiae Mathematicae Japonicae}, Vol. 59 No. 2 (2004). pp. 265-316.

\bibitem{shelah_strong_1}Saharon Shelah. Dependent first order theories, continued. \textit{Israel Journal of Mathematics}, Vol. 173 No. 1 (2009). pp. 1-60.

\bibitem{shelah_strong_2}Saharon Shelah. Strongly dependent theories. \textit{Israel Journal of Mathematics}, Vol. 204 (2014). pp. 1-83.

\bibitem{shelah_g00}Saharon Shelah. Minimal bounded index subgroup for dependent theories. \textit{Proceedings of the AMS}, Vol. 136 No. 3, (2008). pp. 1087-1091.

\bibitem{simon_book}Pierre Simon. A guide to NIP theories. \textit{Cambridge University Press} (2015).

\bibitem{simon_ordered}Pierre Simon. On dp-minimal ordered structures. \textit{The Journal of Symbolic Logic}, Vol. 76 No. 2 (2011). pp. 448-460.

\bibitem{simon_distality}Pierre Simon. Distal and non-distal NIP theories. \textit{Annals of Pure and Applied Logic}, Vol. 164 No. 1 (2013). pp. 294-318.

\bibitem{simonetta_example}Patrick Simonetta. An example of a C-minimal group which is not abelian-by-finite. \textit{Proceedings of the AMS}, Vol. 131 No. 12 (2003). pp. 3913-3917.

\bibitem{simonetta_nilpotent}Patrick Simonetta. On non-abelian C-minimal groups. \textit{Annals of Pure and Applied Logic}, Vol. 122 (2003). pp 263-287.

\bibitem{stonestrom}Atticus Stonestrom. On f-generic types in NIP groups, accepted at \textit{Advances in Mathematics}. arXiv:2303.13470.

\bibitem{usvyatsov_gs}Alexander Usvyatsov. On generically stable types in dependent theories. \textit{The Journal of Symbolic Logic}, Vol. 74, No. 1 (2009). pp. 216-250.

\bibitem{wagner}Frank Wagner. Stable Groups. \textit{Cambridge University Press} (2013).

\end{thebibliography}

\end{document}